\documentclass[12pt]{amsart}

\let\oldlabel=\label
\def\prellabel{\marginparsep=1em\marginparwidth=44pt
    \def\label##1{\oldlabel{##1}\ifmmode\else\ifinner\else
         \marginpar{{\footnotesize\ \\ \tt
                    ##1}}\fi\fi}}
\def\rank{\operatorname{rank}}
\def\depth{\operatorname{depth}}
\def\grade{\operatorname{grade}}
\def\htt{\operatorname{height}}
\def\Ass{\operatorname{Ass}}

\def\gp{\operatorname{gp}}
\def\Ker{\operatorname{Ker}}
\def\Im{\operatorname{Im}}
\def\chara{\operatorname{char}}
\def\Hom{\operatorname{Hom}}

\def\Cl{\operatorname{Cl}}

\def\ini{\operatorname{in}}
\def\supp{\operatorname{supp}}
\def\Supp{\operatorname{Supp}}
\def\Spec{\operatorname{Spec}}
\def\eff{\operatorname{eff}}
\def\Eff{\operatorname{Eff}}

\let\iso\cong
\let\:\colon
\let\Dirsum\bigoplus
\let\dirsum\oplus \let\tensor\otimes

\let\Bbb=\mathbb
\let\frak=\mathfrak
\def\RR{{\Bbb R}}
\def\QQ{{\Bbb Q}}
\def\ZZ{{\Bbb Z}}
\def\NN{{\Bbb N}}
\def\mm{{\frak m}}
\def\pp{{\frak p}}
\def\qq{{\frak q}}
\def\rr{{\frak r}}
\def\FF{{\mathcal F}}
\def\LL{{\mathcal L}}
\def\DD{{\mathcal D}}
\def\DDD{{\overline{\mathcal D}\,'}}
\def\EE{{\mathcal E}}
\def\Cc{{\mathcal C}}
\def\Gc{{\mathcal G}}
\def\Sc{{\mathcal S}}
\def\R{{\mathcal R}}
\def\MM{{\mathcal M}}

\let\epsilon=\varepsilon
\let\phi=\varphi
\let\theta=\vartheta

\let\hat\widehat

\input diagrams
\newarrow{Equal} =====
\newarrow{Into} C--->
\newarrow{Onto} ----{>>}

\newtheorem{lemma}{Lemma}[section]
\newtheorem{corollary}[lemma]{Corollary}
\newtheorem{theorem}[lemma]{Theorem}
\newtheorem{proposition}[lemma]{Proposition}

\theoremstyle{definition}

\newtheorem{remark}[lemma]{Remark}

\newtheorem{example}[lemma]{Example}

\begin{document}
\title[Divisorial linear algebra of normal semigroup rings]
{Divisorial linear algebra\\ of normal semigroup rings}
\author{Winfried Bruns \and Joseph Gubeladze}
\thanks{The second author was partially supported by the Max-Planck-Institut
f\"ur Mathematik at Bonn.}
\address{Universit\"at Osnabr\"uck,
FB Mathematik/Informatik, 49069 Osnabr\"uck, Germany}
\email{Winfried.Bruns@mathematik.uni-osnabrueck.de}
\address{A. Razmadze Mathematical Institute, Alexidze St. 1, 380093
Tbilisi, Georgia} \email{gubel@rmi.acnet.ge}


\begin{abstract}
We investigate the minimal number of generators $\mu$ and the
depth of divisorial ideals over normal semigroup rings. Such
ideals are defined by the inhomogeneous systems of linear
inequalities associated with the support hyperplanes of the
semigroup. The main result is that for every bound $C$ there
exist, up to isomorphism, only finitely many divisorial ideals $I$
such that $\mu(I)\le C$. It follows that there exist only finitely
many Cohen--Macaulay divisor classes. Moreover we determine the
minimal depth of all divisorial ideals and the behaviour of $\mu$
and depth in ``arithmetic progressions'' in the divisor class
group.

The results are generalized to more general systems of linear inequalities
whose homogeneous versions define the semigroup in a not necessarily
irredundant way. The ideals arising this way can also be considered as defined
by the non-negative solutions of an inhomogeneous system of linear diophantine
equations.

We also give a more ring-theoretic approach to the theorem on minimal number of
generators of divisorial ideals: it turns out to be a special instance of a
theorem on the growth of multigraded Hilbert functions.
\end{abstract}

\maketitle
\section{Introduction}

A normal semigroup $S\subset\ZZ^n$ can be described as the set of lattice
points in a finitely generated rational cone. Equivalently, it is the set
\begin{equation*}
S=\{x\in \ZZ^n\: \sigma_i(x)\ge 0,\ i=1,\dots,s\}\tag{$*$}
\end{equation*}
of lattice points satisfying a system of homogeneous inequalities given by
linear forms $\sigma_i$ with integral (or rational) coefficients. For a field
$K$ the $K$-algebra $R=K[S]$ is a normal semigroup ring. In the introduction we
always assume that $S$ is positive, i.e.\ $0$ is the only invertible element in
$S$.

Let $a_1,\dots,a_s$ be integers. Then the set
$$
T=\{x\in \ZZ^n\: \sigma_i(x)\ge a_i,\ i=1,\dots,s\}
$$
satisfies the condition $S+T\subset T$, and therefore the $K$-vector space
$KT\subset K[\ZZ^n]$ is an $R$-module in a natural way.

It is not hard to show that an $R$-module $KT$ is a (fractional) ideal of $R$
if the group $\gp(S)$ generated by $S$ equals $\ZZ^n$. Moreover, if the
presentation $(*)$ of $S$ is irredundant, then the $R$-modules $KT$ are even
divisorial ideals.  These divisorial ideals represent the full divisor class
group $\Cl(R)$. Therefore an irredundant system of homogeneous linear
inequalities is the most interesting from the ring-theoretic point of view, but
we will also treat the general case.

We are mainly interested in two invariants of divisorial ideals $D$, namely
their number of generators $\mu(D)$ and their depth as $R$-modules, and in
particular in the Cohen--Macaulay property. Our main result, based on
combinatorial arguments, is that for each $C\in\ZZ_+$ there exist, up to
isomorphism, only finitely many divisorial ideals $D$ such that $\mu(D)\le C$.
It then follows by Serre's numerical Cohen--Macaulay criterion that only
finitely many divisor classes represent Cohen--Macaulay modules.

The second main result concerns the growth of Hilbert functions of certain
multigraded modules and algebras. Roughly speaking, it says that the Hilbert
function takes values $\le C$ only at finitely many graded components, provided
this holds along each arithmetic progression in the grading semigroup. The
theorem on Hilbert functions can be applied to the minimal number of generators
of divisorial ideals since $R$ can be embedded into a polynomial ring $P$ over
$K$ such that $P$ is a $\Cl(R)$-graded $R$-algebra in a natural way. This leads
to a second proof of the result on number of generators mentioned above.

In Section \ref{Co} we study this ``standard'' embedding of $R$ into a
polynomial ring $P$ over $K$. It has been used many times already, for example
by Hochster \cite{Ho} and Stanley \cite{St1}. The standard embedding into a
polynomial ring graded by the divisor class group, as described in Theorem
2.1(a), are exactly the homogeneous coordinate rings in the sense of [Cox] of
the affine toric varieties $\Spec(k[S])$. Moreover, Theorem 2.1(b) corresponds
to the results in \cite[Section 2]{Cox} on representing (not necessarily
affine) toric varieties as quotients of affine spaces by diagonalizable linear
groups. We include the details because they are essential for the subsequent
sections.

More generally we discuss so-called pure embeddings of normal
affine semigroups into polynomial rings $Q=K[Y_1\dots,Y_s]$, and
the natural splitting of $Q$ into the coset modules of $R$ , i.e.\
those submodules of $Q$ whose monomial basis is defined by a
single residue class modulo $\gp(S)\subset\ZZ^n$. Such pure
embeddings arise from linear actions of diagonalizable linear
algebraic groups on polynomial rings over $K$, and, more
generally, from systems of linear diophantine equations and
congruences. In the setting of invariant theory the coset modules
appear as modules of covariants, and in the context of diophantine
equations they correspond to the sets of non-negative solutions of
the associated inhomogeneous systems. We describe all pure
embeddings of a normal affine semigroup into finitely generated
free semigroups, and characterize those embeddings for which all
coset modules are divisorial ideals.

In Section \ref{TCM} we study divisorial ideals of the form $KT$
where
$$
T=\{z\in\ZZ^n:\sigma_i(z)\ge \sigma_i(\beta)\}
$$
for some $\beta\in \RR^n$. These divisorial ideals, which we call
\emph{conic}, have been shown to be Cohen--Macaulay by Stanley
\cite{St4} and Dong \cite{Do}. The set of conic classes contains
all torsion classes in the divisor class group, but is strictly
larger as soon as $\Cl(R)$ is non-torsion. It turns out that the
conic divisorial ideals are exactly those that appear in the
decomposition of $R$ as a module over its isomorphic image
$R^{(k)}$ under the Frobenius-like endomorphism sending each
monomial to its $k$-th power, $k\in \NN$.

Section \ref{De} gives a combinatorial description of the minimal
depth of all divisorial ideals of $R$: it coincides with the
minimal number of facets $F_1,\dots,F_u$ of the cone generated by
$S$ such that $F_1\cap\dots\cap F_u=\{0\}$.

Section \ref{mu} contains our main result on number of generators. The crucial
point in its proof is that the convex polyhedron $C(D)$ naturally associated
with a divisorial ideal $D$ has a compact face of positive dimension if (and
only if) the class of $D$ is non-torsion. One can show that $\mu(D)\ge
M\lambda$ where $M$ is  a positive constant only depending on the semigroup $S$
and $\lambda$ is the maximal length of a compact $1$-dimensional face of
$C(D)$. Moreover, since the compact $1$-dimensional faces are in discrete
positions and uniquely determine the divisor class, it follows that $\lambda$
has to go to infinity in each infinite family of divisor classes.

The observation on compact faces of positive dimension is also crucial for our
second approach to the number of generators via Hilbert functions. Their
well-established theory allows us to prove quite precise results about the
asymptotic behaviour of $\mu$ and depth along an arithmetic progression in the
divisor class group.

 From the viewpoint of applications it is quite unnatural to study
only the sets of solutions $T(a)$ of inhomogeneous systems of
linear inequalities $\xi_i(x)\ge a_i$, $i=1,\dots,n$, whose
associated homogeneous system defines its set of solutions
$S=T(0)$ irredundantly. Therefore we turn to the general case (but
still under the assumption $\gp(S)=\ZZ^n$) in Section \ref{gen},
replacing divisorial ideals by the \emph{$\xi$-convex ideals}
$I=KT(a)$. The result on numbers of generators has a fully
satisfactory generalization: up to isomorphism over $R$, there
exist only finitely many $\xi$-convex ideals $I$ such that
$\mu(I)\le C$ for every constant $C$.

Section \ref{Gr} finally contains the theorem on the growth of Hilbert
functions outlined above. It is proved by an analysis of homomorphisms of
affine semigroups and their ``modules''.

The Cohen--Macaulay property of coset modules has been
characterized by Stanley \cite{St2, St3} in terms of local
cohomology. Brion \cite{Bri} has shown that the number of
isomorphism classes of Cohen--Macaulay modules of covariants is
finite for certain actions of linear algebraic groups; however,
the hypotheses of his theorem exclude groups with infinitely many
characters. Therefore our result is to some extent complementary
to Brion's.

For unexplained terminology we refer the reader to Bruns and Herzog \cite{BH},
Eisenbud \cite{Ei} or Stanley \cite{St3}. A standard reference for the theory
of divisor class groups and divisorial ideals is Fossum \cite{Fo}. Divisor
class groups of normal semigroup rings are discussed in detail by Gubeladze
\cite{Gu2}. The theory of linear algebraic groups used in Section \ref{Co} can
be found in Humphreys \cite{Hu}.

There also exists a considerable literature on the divisor theory of not
necessarily finitely generated Krull monoids. See Halter-Koch \cite{HK} for a
comprehensive overview and further references.

We are grateful to Cornel Bae\c tica and Robert Koch for their
critical comments that have substantially improved the final
version of this paper.

\section{Coset modules}\label{Co}

Usually we write the operation in a semigroup additively, and in our
terminology a semigroup is always commutative and contains a neutral element
$0$. Let $S\subset T$ be semigroups. Then the \emph{integral closure} (or
\emph{saturation}) of $S$ in $T$ is the subsemigroup
$$
\hat S=\{x\in T\: mx\in S\text{ for some } m\in\ZZ_+,\ m>0\}
$$
of $T$; if $S=\hat S$, then $S$ is \emph{integrally closed} in
$T$. By definition, an \emph{affine semigroup} $S$ is a finitely
generated subsemigroup of $\ZZ^n$. The \emph{normalization}
$\overline S$ of an affine semigroup $S$ is its integral closure
in its group of differences $\gp(S)$. We call an affine semigroup
\emph{normal} if $S=\overline S$. All this terminology is
consistent with (and derived from) its use in commutative algebra:
for a field $K$ the semigroup algebra $K[\hat S]$ is the integral
closure of $K[S]$ in $K[T]$ etc., at least if $T$ can be linearly
ordered and every element in $K[T]$ has a leading monomial.

Let $S\subset \ZZ^n$ be an affine semigroup. One sees easily that
its integral closure in $\ZZ^n$ is
$$
\hat S=C(S)\cap \ZZ^n
$$
where $C(S)\subset \RR^n$, the \emph{cone generated by $S$}, is
the set of all linear combinations of elements of $S$ with
non-negative real coefficients. The normalization $\overline S$ of
$S$ is given by $\bar S=C(S)\cap\gp(S)$. It follows from results
of commutative algebra or Gordan's lemma that $\hat S$ and
$\overline S$ are again affine semigroups.

Since $S$, and therefore the cone $C(S)$, are finitely generated,
$C(S)$ is the intersection of finitely many rational halfspaces,
\begin{equation}
C(S)=\bigcap_{i=1}^s \{x\in\RR^n\: \sigma_i(x)\ge 0\}.\tag{$*$}
\end{equation}
That the halfspaces are rational means that we can define them by
linear forms $\sigma_i\in(\RR^n)^*$ with coprime integer
coefficients.

So far we have considered an arbitrary embedding $S\subset \ZZ^n$.
The smallest lattice of which $S$ is a subsemigroup is $\gp(S)$
which we identify with $\ZZ^r$, $r=\rank S$. If we now require
that the representation $(*)$ (with $n=r$) is minimal, then the
$\sigma_i$ are uniquely determined (up to order). We call them, as
well as their restrictions, to $\ZZ^r$ the {\em support forms} of
$S$, and we set
$$
\supp(S)=\{\sigma_1,\dots,\sigma_s\}.
$$

We will often assume that $S$ is \emph{positive}, i.e.\ $0$ is the only
invertible element in $S$.

Let $K$ be a field and $R=K[S]$. By choosing a basis in $\gp(S)$, we identify
it with $\ZZ^r$ for some $r\ge 0$. From the ring-theoretic point of view, $S$
can then be considered as a semigroup of monomials of the Laurent polynomial
ring $K[X_1^{\pm1},\dots,X_r^{\pm1}]\iso K[\ZZ^r]$. The support forms, so far
defined on $\gp(S)\subset K[X_1^{\pm1},\dots,X_r^{\pm1}]$, can be extended to
discrete valuations on the quotient field of $K[S]$, which we also denote by
$\sigma_i$. The prime ideals
$$
\pp_i=\{x\in K[S]\:\sigma_i(x)\ge 1\}
$$
are exactly the divisorial prime ideals generated by monomials. By a theorem of
Chouinard \cite{Ch}, the divisor class group $\Cl(R)$ is isomorphic to the
quotient of the free abelian group generated by $\pp_1,\dots,\pp_s$ modulo the
subgroup generated by the principal divisors of the monomials $x\in S$. We
define the linear map $\sigma\:\ZZ^r\to\ZZ^s$ by
$$
\sigma(x)=(\sigma_1(x),\dots,\sigma_s(x)).
$$
Note that $\sigma$ is injective if $S$ is \emph{positive}; in fact $x\in S$ and
$-x\in S$ for any element $x\in\Ker\sigma$. One has $\sigma(S)\subset\ZZ_+^s$,
and we can extend $\sigma$ to an embedding
$$
\sigma\: K[S]\to K[\ZZ_+^s]\iso K[Y_1,\dots,Y_s].
$$
We call these embeddings the {\em standard embeddings} of $S$ and $K[S]$
respectively. Note that $\sigma(S)$ (up to an automorphism of $\ZZ_+^s)$ only
depends on the isomorphism class of $S$; in fact, according to Gubeladze
\cite{Gu3}, it only depends on the isomorphism class of the $K$-algebra $K[S]$.

Let $K$ be an algebraically closed field, and $S$ a positive normal affine
semigroup. It is well-known that $K[S]$ can be represented as the ring of
invariants of a diagonal action of a suitable algebraic torus
$T_n(K)=(K^\times)^n$ on a polynomial ring $K[X_1,\dots,X_n]$ for some
$n\in\NN$. Conversely, the ring of invariants of a linear action of a
diagonalizable group $D$ over $K$ on such a polynomial ring has the structure
of a normal affine semigroup ring $K[S]$. The diagonalizable groups are exactly
those isomorphic to a direct product $T_n(K)\times A$ where $A$ is a finite
abelian group whose order is not divisible by $\chara K$. In the following
theorem we will discuss the problem how $D$ and $K[S]$ are related. As we shall
see, the connection is given through the divisor class group of $K[S]$.

\begin{theorem}\label{Standard}
Let $K$ be a field, $S$ a positive normal affine semigroup, $R=K[S]$, and
$\sigma\: R\to P=K[Y_1,\dots,Y_s]$ the standard embedding.
\begin{itemize}
\item[(a)] Then $P$ decomposes as an $R$-module into a direct sum of rank $1$
$R$-modules $M_c$, $c\in \Cl(R)$, such that $M_c$ is isomorphic to
a divisorial ideal of class $c$.
\item[(b)] Suppose that $K$ is algebraically closed and that $\Cl(R)$ does not
contain an element of order $\chara K$. Then $R=P^D$ where
$D=\Hom_\ZZ(\Cl(R),K^\times)\subset T_s(K)$ acts naturally on $P$.
\end{itemize}
\end{theorem}

\begin{proof}
To each divisorial ideal $\pp_1^{(a_1)}\cap\dots\cap\pp_s^{(a_s)}$, $a_i\in
\ZZ$, we associate $(a_1,\dots,a_s)\in\ZZ^s$. Under this assignment, the
principal divisorial ideal generated by $s\in\gp(S)$ is mapped to $\sigma(s)$.
By Chouinard's theorem this yields the isomorphism
$$
\Cl(R)\iso \ZZ^s/\sigma(\ZZ^r).
$$
For each $c\in \ZZ^s/\sigma(\ZZ^r)$ we let $M_c$ be the $K$-vector subspace of
$P$ generated by all monomials whose exponent vector in $\ZZ^s/\sigma(\ZZ^r)$
has residue class $-c$. Then $M_c$ is clearly an $R$-submodule of $P$.
Moreover, by construction, $P$ is the direct sum of these $R$-modules.

It remains to show for (a) that $M_c$, $c\in\Cl(R)$, is a divisorial ideal of
class $c$ (relative to the isomorphism above). We choose a representative
$a=(a_1,\dots,a_s)$ of $c$. Then a monomial $s\in\sigma(\ZZ^r)$ belongs to
$\pp_1^{(a_1)}\cap\dots\cap\pp_s^{(a_s)}$ if and only if $s_i\ge a_i$ for all
$i$, and this is equivalent to $s_i-a_i\in\ZZ_+$ for all $i$. Hence the
assignment $s\mapsto s-a$, in ring-theoretic terms: multiplication by the
monomial $Y^{-a}$, induces an $R$-isomorphism
$\pp_1^{(a_1)}\cap\dots\cap\pp_s^{(a_s)}\iso M_c$.

For (b) we consider the exact sequence
$$
0\to\sigma(\ZZ^r)\to\ZZ^s\to\Cl(R)\to 0.
$$
If $K$ is algebraically closed, then this sequence induces an exact
sequence
$$
1\to \Hom_\ZZ(\Cl(R),K^\times) \to \Hom_\ZZ(\ZZ^s,K^\times) \to
\Hom_\ZZ(\sigma(\ZZ^r),K^\times) \to 1
$$
since $K^\times$ is an injective $\ZZ$-module.

The torus $T_s(K)\iso\Hom_\ZZ(\ZZ^s,K^\times)$ operates naturally on the
Laurent polynomial ring $L=K[Y_1^{\pm1},\dots,Y_s^{\pm1}]$ through the
substitution $Y^t\mapsto \tau(t)Y^t$ for each $\tau\in T_s(K)$ and each
monomial $Y^t$ of $P$, $t\in\ZZ_+^s$. In the following we identify $t$ and
$Y^t$.

The homomorphism
$$
\delta\: \ZZ^s\to \Hom_{\mathrm{alg}}(T_s(K),K^\times),\qquad
\delta(t)(\tau)=\tau(t),
$$
is an isomorphism, and it induces an isomorphism
$$
\Hom_{\mathrm{alg}}\left(\Hom_\ZZ(\sigma(\ZZ^r),K^\times),K^\times\right)\iso U
$$
where $U$ is the subgroup of all $t\in\ZZ^s$ such that $\delta(\tau)(t)=1$ for
all $\tau\in D$. Thus the ring of invariants $L^D$ is the $K$-vector space
spanned by all monomials $t\in U$. There is a commutative diagram
$$
\begin{diagram}
0  & \rTo & \sigma(\ZZ^r) & \rTo & \ZZ^s  & \rTo & \Cl(R) & \rTo & 0\\
   &      & \dTo          &      & \dEqual &      & \dOnto \\
0  & \rTo &  U            & \rTo & \ZZ^s  & \rTo & \Hom_{\mathrm
{alg}}(D,K^\times) & \rTo & 0
\end{diagram}
$$
with exact rows and in which the vertical maps are the natural ones. The
right hand vertical map is an isomorphism if and only if the characteristic of
$K$ does not divide the order of any element of $\Cl(R)$, and therefore
$U=\sigma(\ZZ^r)$ in this case. Then the set of monomials generating $P^D$ is
exactly $\sigma(\ZZ^r)\cap\ZZ_+^s$.
\end{proof}

Let $S$ and $T$ be affine semigroups, $S\subset T$. Then $S$ is called a
\emph{pure subsemigroup} if $S=T\cap\gp(S)$ in $\gp(T)$ (in Hochster's
terminology \cite{Ho} $S$ is a {\em full} subsemigroup). It is not hard to
check that $S$ is a pure subsemigroup if and only if $K[S]\subset K[T]$ is a
pure extension of rings. More generally, we say that an injective homomorphism
$\xi\: S\to T$ is \emph{pure} if $\xi(S)$ is a pure subsemigroup of $T$. If $T$
is normal, then every pure subsemigroup of $T$ is also normal. Examples of pure
embeddings are given by the standard embedding $S\to\ZZ_+^s$, $s=\#\supp(S)$,
and the embedding $S\to\ZZ^n$ resulting from a representation of $K[S]$ as a
ring of invariants of a subgroup of the torus $T_n(K)$ acting diagonally on
$K[Y_1,\dots,Y_n]$.

The most basic source of pure subsemigroups are systems composed of homogeneous
linear diophantine equations
$$
a_{i1}x_1+\dots+a_{in}x_n=0,\qquad i=1,\dots,u,
$$
and homogeneous congruences
$$
b_{j1}x_1+\dots+b_{jn}x_n\equiv 0\mod{c_j}\qquad j=1,\dots,v.
$$
The set of nonnegative solutions $x\in\ZZ_+^n$ of such a system evidently forms
a pure subsemigroup of $\ZZ_+^n$.

In module-theoretic terms the equation $S=T\cap\gp(S)$ says that the monomial
$K$-vector space complement of $K[S]$ in $K[T]$ is a $K[S]$-submodule of
$K[T]$, and this condition is actually equivalent to the requirement that
$K[S]$ is a direct summand of $K[T]$ as a $K[S]$-module.

The next proposition describes all pure embeddings of $S$ into a free affine
semigroup:

\begin{proposition}\label{DirSum}
Let $S$ be a positive normal affine semigroup, and $\xi\: S\to
\ZZ_+^n$ an embedding. Then $\xi$ is a pure embedding if and only
if $n\ge s=\#\supp(S)$ and, up to renumbering, there exist $e_i\in
\ZZ$, $e_i>0$, such that $\xi_i=e_i\sigma_i$, $i=1,\dots,s$.
Moreover, each of $\xi_{s+1},\dots,\xi_n$ is a rational linear
combination of $\sigma_1,\dots,\sigma_s$ with non-negative
coefficients.
\end{proposition}

\begin{proof}
Write $\gp(S)=\ZZ^r$. The embedding $\xi$ can be extended to an $\RR$-linear
embedding $\xi\: \RR^r\to \RR^n$. Then $\xi(C(S))=C(\xi(S))$. Suppose that
$\xi(S)$ is a pure subsemigroup of $\ZZ_+^n$. Then the restrictions of the
coordinate linear forms of $\RR^n$ to $\xi(\RR^r)$ cut out the cone $\xi(C(S))$
from $\xi(\RR^r)$. In other words, the inequalities $\xi_i(s)\ge 0$,
$i=1,\dots,n$, define the cone $C(S)$. Thus the $\xi_i$ belong to the dual cone
$$
C(S)^*=\{\alpha\in\Hom_\RR(\RR^r,\RR)\: \alpha(x)\ge 0\text{ for all }x\in
C(S)\},
$$
and they in fact generate it since $C(S)=C(S)^{**}$. However, then
we must have at least one linear form from each extreme ray of
$C(S)^*$ among the $\xi_i$. The extreme rays of $C(S)^*$
correspond bijectively to the support forms $\sigma_i$, and each
non-zero linear form $\alpha$ in the extreme ray through
$\sigma_i$ is of the form $\alpha=e\sigma_i$ with $e\in\ZZ_+$. The
rest is clear, as well as the converse implication.
\end{proof}

\begin{corollary}\label{SubQuot}
$\Cl(K[S])$ is isomorphic to a subquotient of $\ZZ^n/\xi(\gp(S))$ if $\xi$ is
pure.
\end{corollary}

\begin{proof}
By the proposition we can assume that $\xi_i=e_i \sigma_i$,
$e_i\in\ZZ_+$, $e_i>0$, $i=1,\dots,s$. We define
$\zeta\:\ZZ^n\to\ZZ^s$ as the restriction to the first $s$
coordinates. Then $\theta=\zeta\circ\xi$ is injective, and we have
a chain of subgroups
$$
\theta(\gp(S))\subset \ZZ e_1\dirsum\dots\dirsum\ZZ e_s\subset\ZZ^s.
$$
The corollary follows since $\ZZ^s/\theta(\gp(S))$ is a quotient of
$\ZZ^n/\xi(\gp(S)))$, and $(\ZZ e_1\dirsum\dots\dirsum\ZZ e_s)/\theta(\gp(S))$
is isomorphic to $\Cl(K[S])$.
\end{proof}

Theorem \ref{Standard} and Corollary \ref{SubQuot} show that the divisor class
group sets the lower limit to all diagonalizable groups $D$ having $K[S]$ as
their ring of invariants. In fact, if $D\subset T_n(K)$ acts on
$Q=K[X_1,\dots,X_n]$, then $Q^D$ is generated by a pure (normal) subsemigroup
of the semigroup $\ZZ_+^n$ of monomials of $Q$, and as in the proof of Theorem
\ref{Standard} one sees that $\ZZ^n/\gp(S)\iso\Hom_\ZZ(D,K^\times)$.

If $S$ is a pure subsemigroup of $\ZZ_+^n$, then (as in the special case of the
standard embedding)
$$
Q=K[X_1,\dots,X_n] =\Dirsum_{z\in \ZZ^n/\gp(S)} Q_z
$$
where $Q_z$ is the $K$-vector space generated by all monomials that have
residue class $-z$. Clearly $Q_z$ is an $K[S]$-submodule, and $Q$ is a
$\ZZ^n/\gp(S)$-graded $K[S]$-algebra.

Evidently $Q_z\neq 0$ for all $z\in \ZZ^n/\gp(S)$ if and only if $S$ contains
an element of $\ZZ_+^n$ all of whose components are positive. Every non-zero
$K[S]$-module $Q_z$ is finitely generated and of rank $1$. In fact,
multiplication by the monomial $X^{z}$ maps $Q_z$ into $K[\gp(S)]$ and if the
image contains a monomial with a negative exponent, then $S$ contains an
element whose corresponding component is positive. Since the exponents of the
monomials in $Q_z$ are bounded below by $z$ we can embed $Q_z$ into $K[S]$. The
modules $Q_z$ are called the {\em coset modules} of $S$.

In the invariant-theoretic context discussed above, the $Q_z$ are simply the
modules of covariants associated with the character $z\in
\ZZ^n/\gp(S)\iso\Hom_{\mathrm {alg}}(D,K^\times)$. In fact, an element $x\in
Q=K[X_1,\dots,x_n]$ belongs to $Q_z$ if and only if $\delta(x)=\delta(z)x$ for
all $\delta\in D$.

In the context of systems of linear equations and congruences the coset modules
correspond to inhomogeneous such systems. In fact, the monomials forming the
$K$-basis of $Q_z$ are the set of non-negative solutions of an inhomogeneous
system of linear diophantine equations and inequalities whose associated
homogeneous system defines $S$.

\begin{proposition}\label{divi}
Suppose that $S$ is a positive normal affine semigroup and $\xi\: S\to \ZZ_+^n$
a pure embedding. Let $R=K[S]$, $s=\#\supp(S)$, and $Q=K[X_1,\dots,X_n]$. With
the notation as in Proposition \ref{DirSum}, the following are equivalent:
\begin{itemize}
\item[(a)] every non-zero coset module $Q_z$ of $\xi(S)$ is divisorial;
\item[(b)] for each $i=1,\dots,n$ there exists $j_i\in \{1,\dots,s\}$ and
$e_i\in\ZZ_+$ with $\xi_i=e_i\sigma_{j_i}$.
\end{itemize}
\end{proposition}

\begin{proof}
The implication (b)$\implies$(a) is proved in the same way as Theorem
\ref{Standard}. In fact,
$$
Q_z\iso \bigcap_{i=1}^n \pp_{j_i}^{(z_i')},\qquad z_i'=\lceil z_i/e_i\rceil.
$$

For the converse implication we may assume that $\xi_n$ is not a multiple of
any $\sigma_j$. Then we write $\xi_n$ as a non-negative linear combination of
$\sigma_1,\dots,\sigma_s$ with as few non-zero coefficients as possible. By
Carath\'edory's theorem on convex cones and after reordering the $\sigma_j$ we
can assume that $\xi_n=a_1\sigma_1+\dots+a_p\sigma_p$ where $p\ge 2$,
$\sigma_1,\dots,\sigma_p$ are linearly independent, and $a_1,\dots,a_p>0$ are
rational numbers.

It is not hard to check that $\gp(S)$ contains elements $t_j$ with
$\sigma_j(t_j)<0$ and $\sigma_k(t_j)\ge 0$ for all $k=1,\dots,s$, $j\neq k$.
Furthermore $S$ contains elements $u_1,\dots,u_s$ with $\sigma_j(u_j)=0$ and
$\sigma_k(u_j)>0$ for $j\neq k$. Given $q\in\ZZ_+$ and $j\in\{1,\dots,s\}$, we
can therefore find an element $v_{jq}\in\gp(S)$ such that $\sigma_j(v_{jq})\le
-q$, but $\xi_n(v_{jq})\ge0$.

Set
$$
z_m=(-m,\dots,-m,0) \mod \xi(\gp(S)).
$$
Then $Q_{z_m}$ is isomorphic to the fractional ideal $I(m)$ spanned by all the
monomials $x\in S$ with $\xi_j(x)\ge -m$ for $j=1,\dots,n-1$ and $\xi_n(x)\ge
0$. Let
$$
J(m)=\pp_1^{(b_1(m))}\cap\dots\cap \pp_s^{(b_s(m))}
$$
be the smallest divisorial ideal containing $I(m)$. Then the existence of the
elements $v_{jq}$ implies $b_i(m)\to-\infty$ with $m\to\infty$ for all
$i=1,\dots,s$.

On the other hand we can find an element $w\in\gp(S)$ with $\sigma_1(w)<0$ and
$\sigma_j(w)=0$ for $j=2,\dots,p$. In particular $w\in J(m)$ for $m\gg0$.
However $\xi_n(w)<0$. Therefore $I(m)\neq J(m)$ for $m\gg0$, and $Q_{z_m}\iso
I(m)$ is not divisorial.
\end{proof}

Condition (b) of Proposition \ref{divi} can easily be checked algorithmically.
For example, if a system of generators $E=\{x_1,\dots,x_m\}$ of $S$ is known,
then one forms the sets $E_i=\{x_j\: \xi_i(x_j)=0\}$ for each coordinate
$i=1,\dots, n$. Condition (b) is satisfied if and only if none of the $E_i$
contains some $E_k$ properly, and one has $n=\#\supp(S)$ if and only if none of
the $E_i$ contains some $E_k$, $k\neq i$.

\section{Divisor classes associated with torsion cosets}\label{TCM}

The following example illustrates several theorems proved in this
and the following sections.

\begin{example}\label{det2}
Consider the Segre product
$$
R_{mn}=K[X_iY_j\: 1\le i\le m,\ 1\le j \le n]\subset
P=K[X_1,\dots,X_m,Y_1,\dots,Y_n]
$$
of the polynomial rings $K[X_1,\dots,X_m]$, $m\ge 2$, and
$K[Y_1,\dots,Y_n]$, $n\ge 2$, with its standard embedding. It has
divisor class group isomorphic to $\ZZ$, and the two generators of
$\Cl(R)$ correspond to the coset modules $M_{1}=RX_1+\dots+RX_m$
and $M_{-1}=RY_1+\dots+RY_n$. Therefore
$$
\mu(M_{i})=\binom{m+i-1}{m-1},\qquad \mu(M_{-i})=\binom{n+i-1}{n-1}.
$$
for all $i\ge 0$. The Cohen--Macaulay divisorial ideals are
represented by
$$
M_{-(m-1)},\dots,M_0=R_{mn},\dots,M_{n-1}
$$
(see Bruns and Guerrieri \cite{BGue}), and in particular, their
number is finite. However, the finiteness of the number of
Cohen--Macaulay classes is not a peculiar property of $R_{mn}$: it
holds for all normal semigroup rings, as we will see in Corollary
\ref{CMfin}.

Moreover, one has
$$
\inf_{i\ge 0} \depth M_{i}=n,\qquad \inf_{i\ge 0} \depth M_{-i}=m,
$$
and for $i\gg0$ the minimal values are attained (see
\cite[(9.27)]{BV}). Set $p(i)=\mu(M_i)$. It follows that the
degree of the polynomial $p$ and $\inf_i \depth M_i$ add up to
$m+n-1=\dim R$. This is another instance of a general fact (see
Theorem \ref{aripro}).
\end{example}

As a positive result on the Cohen--Macaulay property we now prove
that cosets of torsion elements relative to pure embeddings yield
Cohen--Macaulay divisorial ideals. The method of proof has been
used before by Stanley for the derivation of an analogous result
on semi-invariants of torus actions \cite[3.5]{St4}. Although the
following theorem is essentially equivalent to Stanley's result,
it might be useful to include a discussion of the Cohen--Macaulay
divisorial ideals arising from it.

\begin{theorem}\label{torsion}
Let $S\subset A$ be a pure extension of affine normal semigroups
and set $R=K[S]$, $Q=K[A]$. Furthermore let $\hat S=\{z\in A\:
mz\in S\text{ for some }m\in \NN\}$ be the integral closure of $S$
in $A$.
\begin{itemize}
\item[(a)] For all $z\in \gp(\hat S)/\gp(S)$ the coset module $Q_z$ is a
Cohen--Macaulay $R$-module, provided $Q_z\neq 0$.
\item[(b)] In particular every divisorial ideal $I$ of $R$ whose class in
$\Cl(R)$ is a torsion element is a Cohen--Macaulay $R$-module.
\end{itemize}
\end{theorem}

\begin{proof}
(a) Since $\hat S$ is a normal affine semigroup, the ring $K[\hat
S]$ is Cohen--Macaulay by Hochster's theorem \cite{Ho}. It
decomposes into the direct sum of the finitely many and finitely
generated $R$-modules $Q_z$ considered in (a). Therefore it is a
Cohen--Macaulay ring if and only if all the $Q_z$ are
Cohen--Macaulay $R$-modules.

(b) This follows from (a) if one chooses $S\subset\ZZ_+^s$ as the standard
embedding of $R$.
\end{proof}

The introductory example shows that also non-torsion classes may
be Cohen--Macaulay, and as we will see now, the divisorial ideals
covered by part (a) of Theorem \ref{torsion} may very well have
non-torsion classes in $\Cl(R)$. The following corollary was
proved by Dong \cite{Do} using topological methods:

\begin{corollary}\label{Dong}
Let $S\subset \ZZ^n=\gp(S)$ be a normal positive affine semigroup,
and set and $R=K[S]$. Then the $R$-module $KT$ is a divisorial
Cohen--Macaulay ideal for all subsets $T=\ZZ^n\cap(\beta+C(S))$ of
$\ZZ^n$, $\beta\in\RR^n$.
\end{corollary}

\begin{proof}
Set $a_i=\lceil \sigma_i(\beta)\rceil$ for $i=1,\dots,s$ (where,
as usual $\sigma_1,\dots,\sigma_s$ are the support forms of $S$).
Then clearly
$$
KT=\pp_1^{(a_1)}\cap\dots\cap \pp_s^{(a_s)},
$$
and so $KT$ is a divisorial ideal. It is clear that we can assume
$\beta\in\QQ^n$. Suppose that $m\beta\in\ZZ^n$. Then we consider
the pure embedding $m\sigma$ of $R$ into $K[Y_1,\dots,Y_s]$. The
assignment
$$
z\mapsto m\sigma(z)-\sigma(m\beta)
$$
maps $T$ to $\bigl(m\sigma(\ZZ^n)-\sigma(m\beta)\bigr)\cap
\ZZ_+^s$, and therefore establishes an $R$-isomorphism of $KT$
with the coset module associated with $\sigma(m\beta)$.  But
$m\sigma(m\beta)\in\sigma(\ZZ^n)$, and so part (a) of the theorem
applies.

For use below we notice that the coset $z$ for which the
associated module is isomorphic to $KT$ has been realized within
$\sigma(\ZZ^n)$.
\end{proof}

We call the divisorial ideals and divisor classes considered in
the corollary \emph{conic}.

\begin{remark}
(a) Among the conic divisorial ideals are the ideals $\qq_G$ defined by a face
$G$ of $C(S)$ in the following way:
$$
\qq_G=\bigcap_{F_i\supset G} \pp_i.
$$
This includes the divisorial prime ideals $\pp_1,\dots,\pp_s$ and their intersection
$\omega=\pp_1\cap\dots\cap\pp_s$, which is well-known to be the canonical module of $R$.
(Of course, these are known to be Cohen--Macaulay for other reasons.) Therefore,
whenever $\Cl(R)$ is not a torsion group, then there exist conic
divisor classes that have non-torsion class.

In order to prove that $\qq_G$ is conic we have to find
$\beta\in\RR^n$ such that $0<\sigma_i(\beta)\le 1$ if $F_i\supset
G$ and $-1<\sigma_i(\beta)\le 0$ otherwise. Choose $\gamma$ in the
interior of $-C(S)$ and $\delta$ in the relative interior of $G$.
Then $\sigma_i(\gamma+c\delta)$ has the desired sign for $c\gg0$
and we choose $\beta=c'(\gamma+c\delta)$ for sufficiently small
$c'>0$.

(b) For every Cohen--Macaulay divisorial ideal $\qq$ the ideal
$\omega:\qq=\Hom_R(\qq,\omega)$ is also Cohen--Macaulay (see \cite[3.3.10]{BH}). It is
interesting to observe that $\omega:\qq$ is conic if $\qq$ is so. In fact, suppose that
$\qq$ corresponds to $\ZZ^n\cap(\beta+C(S))$. Then we can assume that
$\sigma_i(\beta)\notin \ZZ$, and, replacing $\beta$ by $\beta+\gamma$ for some interior
monomial $\gamma$ of $S$, we can also assume that $\sigma_i(\beta)>0$, without changing
the class of $\qq$. Then $\lceil
\sigma_i(-\beta)\rceil=1-\lceil\sigma_i(\beta)\rceil$, and thus $\omega:\qq$ is isomorphic
to the conic divisorial ideal determined by $-\beta$, and therefore conic itself.

 From (a) we conclude that the divisorial ideals
$$
\rr_G=\bigcap_{F_i\not\supset G} \pp_i
$$
are also conic.
\end{remark}

As Dong noticed, the previous corollary immediately implies the
following theorem of Stanley (see \cite[3.5]{St4} and
\cite[3.2]{St2}):

\begin{corollary}\label{stan}
Let $\Phi:\ZZ^s\to\ZZ^u$ be a $\ZZ$-linear map. Let $K$ be a
field, and let $R$ be the semigroup ring defined by the
non-negative solutions of the system $\Phi(y)=0$. Suppose there is
a real solution $\beta$ of the system $\Phi(\beta)=\alpha$ for
$\alpha\in\ZZû$ such that
\begin{itemize}
\item[(i)] $-1<\beta_i\le 0$ for all $i$,
\item[(ii)] $z_i\ge0$ for $i=1\dots,s$ if $z$ is an integral solution of
$\Phi(z)=\alpha$ with $z_i\ge \beta_i$ for $i=1,\dots,s$.
\end{itemize}
Then the set of non-negative integral solutions of the system
$\Phi(z)=\alpha$ defines a Cohen--Macaulay divisorial $R$-module.
\end{corollary}

For each semigroup ring $R=K[S]$ there exist natural endomorphisms
$\iota^{(k)}:R\to R$, given by the assignment $s\mapsto s^k$ for
all $s\in S$ (in multiplicative notation). For affine semigroups
$\iota^{(k)}$ is evidently an embedding, and we denote the
(isomorphic) image of $R$ by $R^{(k)}$. Then $R$ is a finitely
generated $R^{(k)}$-module. The pure extensions $\iota^{(k)}$
yield all the divisor classes discussed above:

\begin{proposition}\label{join}
Let $S$ and $R$ be as in Corollary \ref{Dong}. Then the following
constructions all lead to the same set of (Cohen--Macaulay)
divisor classes:
\begin{itemize}
\item[(a)] Theorem \ref{torsion}(a),
\item[(b)] the conic classes,
\item[(c)] the set of divisorial ideals arising from the decomposition
of $R$ as an $R^{(k)}$-module, $k\in \NN$,
\item[(d)] the set of divisorial ideals arising from the decomposition
of $R$ as an $R^{(k)}$-module for a single $k$, provided $k\gg 0$.
\end{itemize}
\end{proposition}

\begin{proof}
Class (c) is obviously contained in class (a), and we proved
Corollary \ref{Dong} by showing that class (b) is contained in
class (c). With the methods we have discussed already, it is easy
to see that class (a) is contained in (b), and so all three
classes (a), (b), (c) coincide.

Up to isomorphism, each conic class can be realized as $KT$ with
$T=\ZZ^n\cap (\beta+C(S))$ with $-1<\beta_i\le0$ since
translations by integral vectors preserve the isomorphism class.
Since the number of possible values $\lceil\sigma_i(\beta)\rceil$
is finite for the $\beta$ under consideration, there exist only
finitely many conic classes. (This results also from the
finiteness of the number of Cohen--Macaulay divisor classes.) Now
it follows easily that suitable vertices $\beta$ can be chosen
simultaneously in $\ZZ^n/k$, provided $k$ is sufficiently large.
\end{proof}

It is an interesting observation that for the example $R=R_{mn}$
discussed in \ref{det2} all Cohen--Macaulay classes are conic: for
$i$ and $k$ let $N_{ki}$ the $R^{(k)}$- submodule of
$P=K[X_1,\dots,X_m,Y_1,\dots,Y_n]$ generated by the $i$-th powers
of the degree $k$ monomials in $X_1,\dots,X_m$. As an $R\iso
R^{(k)}$-module $N_{ki}$ is isomorphic to $M_i$. Now consider a
monomial $\mu=Y_1^{e_1}\cdots Y_n^{e_n}$  of degree $ki$ in
$Y_1,\dots,Y_n$. If $i\le n-1$ and $k$ is sufficiently big, then
we can choose $e_1,\dots,e_n\le k-1$, and it is not hard to see
that $\mu N_{ki}\iso N_{ki}\iso M_i$ is a coset module of the pure
extension $R^{(k)}\to R$. Similarly one realizes the $M_{-i}$ for
$0\le i\le m-1$.

\section{Asymptotic depth of divisor classes}\label{De}

The next theorem describes the asymptotic behaviour of the depth
of divisorial ideals. (See the proof of Theorem \ref{Standard} for
the definition of the modules $M_c$ representing the divisor
classes.) By $\depth_R M$ or simply $\depth M$ we denote the
length of a maximal $M$-sequence for a finitely generated
$R$-module $M$. If $S$ is positive, $R=K[S]$ and $M$ is
$\gp(S)$-graded, then $\depth_R M=\depth_{R_\mm} M_\mm$ where
$\mm$ is the maximal ideal generated by all the non-unit monomials
of $R$. (See \cite[1.5.15]{BH} for the proof of an analogous
statement about $\ZZ$-graded rings and modules.) For an ideal $I$
in a noetherian ring $R$ we let $\grade I$ denote the length of a
maximal $R$-sequence in $I$; one always has $\grade I\le\htt I$.
But if $R$ is Cohen--Macaulay, then $\grade I=\htt I$.

\begin{theorem}\label{depth}
Let $K$ be a field, $S$ a positive normal affine semigroup, $R=K[S]$, and
$\sigma\: R\to P=K[X_1,\dots,X_s]$ the standard embedding. Furthermore let
$\mm$ be the irrelevant maximal ideal of $R$ generated by all non-unit
monomials, and $\lambda$ the maximal length of a monomial $R$-sequence. Then
$$
\lambda\le \grade \mm P = \min\{\depth M_c\: c\in\Cl(R)\}.
$$
\end{theorem}

\begin{proof}
For the inequality it is enough to show that a monomial $R$-sequence is also a
$P$-sequence. (It is irrelevant whether we consider $P$ as an $R$-module or a
$P$-module if elements from $R$ are concerned.) Let $\mu_1,\dots,\mu_u$ be
monomials in $R$ forming an $R$-sequence. Then the subsets
$A_i=\Ass_R(R/(\mu_i))$ are certainly pairwise disjoint. On the other hand,
$A_i$ consists only of monomial prime ideals of height $1$ in $R$, since $R$ is
normal. So $A_i=\{\pp_j\: \sigma_j(\mu_i)>0\}$, and the sets of indeterminates
of $P$ that divide $\sigma(\mu_i)$ in $P$, $i=1,\dots,u$, are pairwise
disjoint. It follows that $\mu_1,\dots,\mu_u$ form a $P$-sequence.

In order to prove the equality we first extend the field $K$ to an uncountable
one. This is harmless, since all data are preserved by base field extension.
Then we can form a maximal $P$-sequence in $\mm P$ by elements from the
$K$-vector subspace $\mm$. Such a $P$-sequence of elements in $R$ then has
length equal to $\grade \mm P$ and is clearly an $M$-sequence for every
$R$-direct summand $M$ of $P$, and in particular for each of the modules $M_c$
representing the divisor classes. Thus $\depth M_c\ge \grade \mm P$.

Whereas this argument needs only finite prime avoidance, we have to use
countable prime avoidance for the converse inequality. Suppose that
$u<\min\{\depth M_c\: c\in\Cl(R)\}$ and that $x_1,\dots,x_u$ is a $P$-sequence
in $\mm$. Then the set
$$
A=\bigcup_{c\in\Cl(R)}\Ass\bigl(M_c/(x_1,\dots,x_u)M_c\bigr)
$$
is a countable set of $K$-vector subspaces of $\mm$. Each prime ideal
associated to $M_c/(x_1,\dots,x_u)M_c$ is a proper subspace of $\mm$ because of
$u<\depth M_c$. Hence $A$ cannot exhaust $\mm$, as follows from elementary
arguments. So we can choose an element $x_{u+1}$ in $\mm$ not contained in a
prime ideal associated to any of the $M_c/(x_1,\dots,x_u)M_c$. So $x_{u+1}$
extends $x_1,\dots,x_u$ to an $M_c$-sequence simultaneously for all
$c\in\Cl(R)$.
\end{proof}

Both the numbers $\lambda$ and $\grade \mm P$ can be characterized
combinatorially:

\begin{proposition}\label{combi}
With the notation of the previous theorem, the following hold:
\begin{itemize}
\item[(a)]
$\grade \mm P$ is the minimal number $u$ of facets $F_{i_1},\dots,F_{i_u}$ of
$C(S)$ such that $F_{i_1}\cap\dots\cap F_{i_u}=\{0\}$.
\item[(b)]
$\lambda$ is the maximal number $\ell$ of subsets $\FF_1,\dots,\FF_\ell$ of
$\FF=\{F_1,\dots,F_s\}$ with the following properties:
$$
\mathrm{(i)}\quad \FF_i\cap\FF_j=\emptyset,\qquad \mathrm{(ii)}
\bigcap_{F\in\FF\setminus \FF_i} F\neq\{0\}
$$
for all $i,j$ such that $i\neq j$.
\end{itemize}
\end{proposition}

\begin{proof}
(a) The ideal $\mm P$ of $P$ is generated by monomials. Therefore all its
minimal prime ideals are generated by indeterminates of $P$. The ideal
generated by $X_{i_1},\dots,X_{i_u}$ contains $\mm P$ if and only if for each
monomial $\mu\in\mm$ there exists a $\sigma_{i_j}$ such that
$\sigma_{i_j}(\mu)>0$. The monomials for which none such inequality holds are
precisely those in $F_{i_1}\cap\dots\cap F_{i_u}$.

(b) Let $\mu_1,\dots,\mu_\ell$ be a monomial $R$-sequence. Then the sets
$\FF_i=\{F\:\pp_F\in\Ass(R/(\mu_i))\}$ are pairwise disjoint, and moreover
$\mu_i\in \bigcap_{F\in\FF\setminus \FF_i} F$. Thus conditions (i) and (ii) are
both satisfied.

For the converse one chooses monomials $\mu_i\in\bigcap_{F\in\FF\setminus
\FF_i} F$. Then $\{F\:\pp_F\in\Ass(R/(\mu_i))\}\subset\FF_i$, and since the
$\FF_i$ are pairwise disjoint, the $\mu_i$ form even a $P$-sequence as observed
above.
\end{proof}

\begin{remark}\label{simpl}
The following are equivalent:
\begin{itemize}
\item[(a)] $S$ is simplicial;
\item[(b)] every divisorial ideal of $R$ is a Cohen--Macaulay module.
\end{itemize}

If $S$ is simplicial, then $\Cl(R)$ is finite, and so every divisorial ideal is
Cohen--Macaulay by Theorem \ref{torsion}. Conversely, if every divisorial ideal
is Cohen--Macaulay, then it is impossible that facets $F_1,\dots,F_u$, $u<\rank
S$, of $C(S)$ intersect only in $0$, and this property characterizes simplicial
cones. However, the implication (b)$\implies$(a) follows also from Theorem
\ref{finite} below.

It is an amusing consequence that a rank $2$ affine semigroup is simplicial --
an evident geometric fact. In fact, if $\dim R\ge 2$, then $\depth I\ge 2$ for
every divisorial ideal.
\end{remark}

The results of this section can be partially generalized to arbitrary pure
embeddings; we leave the details to the reader.

\section{The number of generators}\label{mu}

In this section we first prove our main result on the number of generators of
divisorial ideals of normal semigroup rings $R=K[S]$. In its second part we
then show that it can also be understood and proved as an assertion about the
growth of the Hilbert function of a certain multigraded $K$-algebra.

\begin{theorem}\label{finite}
Let $R$ be a positive normal semigroup ring over the field $K$, and
$m\in\ZZ_+$. Then there exist only finitely many $c\in\Cl(R)$ such that a
divisorial ideal $D$ of class $c$ has $\mu(D)\le m$.
\end{theorem}

As a consequence of Theorem \ref{finite}, the number of Cohen--Macaulay classes
is also finite:

\begin{corollary}\label{CMfin}
There exist only finitely many $c\in \Cl(R)$ for which a divisorial ideal of
class $c$ is a Cohen--Macaulay module.
\end{corollary}

\begin{remark}\label{finrem}
(a) One should note that $\mu(D)$ is a purely combinatorial invariant. If $S$
is the underlying semigroup and $T$ is a monomial basis of a monomial
representative of $D$, then $\mu(D)$ is the smallest number $g$ such that there
exists $x_1,\dots,x_g\in T$ with $T=(S+x_1)\cup\dots\cup(S+x_g)$. Therefore
Theorem \ref{finite} can very well be interpreted as a result on the generation
of the sets of solutions to inhomogeneous linear diophantine equations and
congruences (with fixed associated homogeneous system).

(b) Both the theorem and the corollary hold for all normal affine
semigroups $S$, and not only for positive ones. It is easily seen
that a normal semigroup $S$ splits into a direct summand of its
largest subgroup $S_0$ and a positive normal semigroup $S'$. Thus
one can write $R=K[S]$ as a Laurent polynomial extension of the
$K$-algebra $R'=K[S']$. Each divisor class of $R$ has a
representative $D'\tensor_{R'} R$. Furthermore
$\mu_{R'}(D')=\mu_R(D'\tensor_{R'} R)$ and the Cohen--Macaulay
property is invariant under Laurent polynomial extensions.
\end{remark}

We first derive the corollary from the theorem. Let $\mm$ be the irrelevant
maximal  ideal of $R$. If $M_c$ is a Cohen--Macaulay module, then $(M_c)_\mm$
is a Cohen--Macaulay module (and conversely). Furthermore
$$
e((M_c)_\mm)\ge \mu((M_c)_\mm)=\mu(M_c).
$$
By Serre's numerical Cohen--Macaulay criterion (for example, see
\cite[4.7.11]{BH}), the rank $1$ $R_\mm$-module $(M_c)_\mm$ is Cohen--Macaulay
if and only if its multiplicity $e((M_c)_\mm)$ coincides with $e(R_\mm)$.

\begin{proof}[Proof of Theorem \ref{finite}]
Let $D$ be a monomial divisorial ideal of $R$. As pointed out already, there
exist integers $a_1,\dots,a_s$ such that the lattice points in the set
$$
C(D)=\{x\in \RR^r\:\sigma_i(x)\ge a_i,\ i=1,\dots,s\}
$$
give a $K$-basis of $D$ (again $s=\#\supp(S)$, and the $\sigma_i$
are the support forms). The polyhedron $C(D)$ is uniquely
determined by its extreme points since each of its facets is
parallel to one of the facets of $C(S)$ and passes through such an
extreme point. (Otherwise $C(D)$ would contain a full line, and
this is impossible if $S$ is positive.)

Moreover, $D$ is of torsion class if and only if $C(D)$ has a single extreme
point. This has been proved in \cite{Gu2}, but since it is the crucial point
(sic!) we include the argument.

Suppose first that $D$ is of torsion class. Then there exists $m\in\ZZ_+$,
$m>0$, such that $D^{(m)}$ is a principal ideal, $D^{(m)}=xR$ with a monomial
$x$. It follows that $C(D^{(m)})=mC(D)$ has a single extreme point in (the
lattice point corresponding to) $x$, and therefore $C(D)$ has a single extreme
point.

Conversely, suppose that $C(D)$ has a single extreme point. The extreme point
has rational coordinates. After multiplication with a suitable $m\in\ZZ_+$,
$m>0$, we obtain that $C(D^{(m)})=mC(D)$ has a single extreme point $x$ which
is even a lattice point. All the facets of $C(D^{(m)})$ are parallel to those
of $S$ and must pass through the single extreme point. Therefore $C(D^{(m)})$
has the same facets as $C(S)+x$. Hence $C(D^{(m)})=C(S)+x$. This implies
$D^{(m)}=Rx$ (in multiplicative notation), and so $m$ annihilates the divisor
class of $D$.

Suppose that $D$ is not of torsion class. We form the line complex
$\overline\LL$ consisting of all $1$-dimensional faces of the polyhedron
$C(D)$. Then $\overline\LL$ is connected, and each extreme point is an endpoint
of a $1$-dimensional face. Since there are more than one extreme points, all
extreme points are endpoints of compact $1$-dimensional faces, and the line
complex $\LL(D)$ formed by the \emph{compact} $1$-dimensional faces is also
connected. Since each facet passes through an extreme point, $D$ is uniquely
determined by $\LL(D)$ (as a subset of $\RR^s$).

Let $\Cc$ be an infinite family of divisor classes and choose a divisorial
ideal $D_c$ of class $c$ for each $c\in \Cc$. Assume that the minimal number of
generators $\mu(D_c)$, $c\in \Cc$, is bounded above by a constant $C$. By Lemma
\ref{length} below the Euclidean length of all the line segments
$\ell\in\LL(D_c)$, $c\in \Cc$, is then bounded by a constant $C'$.

It is now crucial to observe that the endpoints of all the line segments under
consideration lie in an overlattice $L=\ZZ^n[1/d]$ of $\ZZ^n$. In fact each
such point is the unique solution of a certain system of linear equations
composed of equations $\sigma_i(x)=a_i$, and therefore can be solved over
$\ZZ[1/d]$ where $d\in\ZZ$ is a suitable common denominator. (Again we have
denoted the support forms of $S$ by $\sigma_i$.)

Let us consider two line segments $\ell$ and $\ell'$ in $\RR^n$ as equivalent
if there exists $z\in \ZZ^n$ such that $\ell'=\ell+z$. Since the length of all
the line segments under consideration is bounded and their endpoints lie in
$\ZZ^n[1/d]$, there are only finitely many equivalence classes of line segments
$\ell\in\LL(D_c)$, $c\in\Cc$.

Similarly we consider two line complexes $\LL(D)$ and $\LL(D')$ as equivalent
if $\LL(D')=\LL(D)+z$, $z\in\ZZ^n$. However, this equation holds if and only if
$C(D')=C(D)+z$, or, in other words, the divisor classes of $D$ and $D'$
coincide.

Since there are only finitely many equivalence classes of line segments and the
number of lines that can appear in a complex $\LL(D)$ is globally bounded (for
example, by $2^s$), one can only construct finitely many connected line
complexes that appear as $\LL(D)$, up to equivalence of line complexes. This
contradicts the infinity of the family $\Cc$.
\end{proof}

\begin{lemma}\label{length}
Let $S$ be a positive normal semigroup, $K$ a field, $D$ a monomial divisorial
ideal whose class is not torsion. Then there exists a constant $M>0$, which
only depends on $S$, such that $\mu(D)\ge M\lambda$ where $\lambda$ is the
maximal Euclidean length of a compact $1$-dimensional face of $C(D)$.
\end{lemma}

\begin{proof}
We assume that $\ZZ^n=\gp(S)$ so that the cone $C(S)$ and the polyhedron $C(D)$
are subsets of $\RR^n$.  Let $\ell$ be a $1$-dimensional compact face of
$C(D)$. Suppose $D$ is given by the inequalities
$$
\sigma_i(x)\ge a_i,\qquad i=1,\dots,s\qquad \bigl(s=\#\supp(C(S))\bigr).
$$
There exists $\epsilon>0$ such that $U_\epsilon(x)\cap C(D)$ contains a lattice
point for each $x\in C(D)$. (In fact, $C(S)$ contains a unit cube, and
$x+C(S)\subset C(D)$ for $x\in C(D)$.) Let $x\in\ell$. We can assume that
$$
\sigma_i(x)\begin{cases}=a_i,&i=1,\dots,m,\\
                         >a_i,&i>m.
            \end{cases}
$$
Let $\tau=\sigma_1+\dots+\sigma_m$. There exists $C>0$ such that $\tau(y)<C$
for all $y\in\RR^n$ with $|y|<\epsilon$.

Furthermore we have $\tau(z)>0$ for all $z\in C(S)$, $z\neq 0$. Otherwise the
facets $F_1,\dots,F_m$ would meet in a line contained in $C(S)$, and this is
impossible if $\ell$ is compact. In particular there are only finitely many
lattice points $z$ in $S$ such that $\tau(z)<C$, and so there exists $\delta>0$
such that $\tau(z)<C$ for $z\in S$ is only possible with $|z|<\delta$.

Now suppose that $D$ is generated by $x_1,\dots,x_q$. For $x\in\ell$ we choose
a lattice point $p\in U_\epsilon(x)\cap C(D)$.  By assumption there exists
$z\in C(S)$ such that $p=x_i+z$. Then
$$
\tau(z)=\tau(p)-\tau(x_i)\le \tau(p)-\tau(x)=\tau(p-x)< C.
$$
Thus $|z|<\delta$, and therefore $|x-x_i|<\delta+\epsilon$.

It follows that the Euclidean length of $\ell$ is bounded by
$2q(\delta+\epsilon)$. Of course $\delta$ depends on $\tau$, but there exist
only finitely many choices for $\tau$ if one varies $\ell$.
\end{proof}

As pointed out, the polyhedron $C(D)$ contains a $1$-dimensional compact face
if $D$ is not of torsion class, but in general one cannot expect anything
stronger. On the other hand, there exist examples for which $C(D)$ for every
non-torsion $D$ has a compact face of arbitrarily high dimension; see Example
\ref{det2}.

If $C(D)$ has a $d$-dimensional face $F$, then the argument in the proof of
Lemma \ref{length} immediately yields that $\mu(D^{(j)})\ge Mj^d$ for a
constant $M>0$: one has only to replace the length of the line segment by the
$d$-dimensional volume of $F$. We now give another proof of a slightly more
general statement. As we will see, it leads to a quite different proof of
Theorem \ref{finite}.

Let $S$ be a positive normal affine semigroup. Recall that the polynomial ring
$P$ of the standard embedding $\sigma\: R\to P$ decomposes into the direct sum
of modules $M_c$, $c\in\Cl(R)$. In the following $C(M_c)$ stands for any of the
congruent polyhedra $C(D)$ where $D$ is a divisorial ideal of class $c$.

\begin{theorem} \label{aripro}
Let $c,d\in\Cl(R)$ and suppose that $c$ is not a torsion element.
\begin{itemize}
\item[(a)]
Then $\lim_{j\to\infty}\mu(M_{jc+d})=\infty$.
\item[(b)]
More precisely, let $m$ be the maximal dimension of the compact faces of
$C(M_c)$. Then there exists $e\in\NN$ such that
$$
\lim_{j\to\infty}\mu(M_{(ej+k)c+d})\frac{m!}{j^{m}}
$$
is a positive natural number for each $k=0,\dots,e-1$.
\item[(c)] One has $\inf_j\depth M_{jc}=\dim R-m$ and $\inf_j\depth
M_{cj+d}\le\dim R-m$.
\end{itemize}
\end{theorem}

\begin{proof}
Let
$$
\DD=\Dirsum_{j=0}^\infty M_{jc}\qquad\text{and}\qquad \MM=\Dirsum_{j=0}^\infty
M_{jc+d}.
$$
Then $\DD$ is a finitely generated $K$-algebra. This follows for general
reasons from Theorem \ref{fingen} below: $\DD$ is the direct sum of graded
components of the $\Cl(R)$-graded $R$-algebra $P$, taken over a finitely
generated subsemigroup of $\Cl(R)$. Theorem \ref{fingen} also shows that $\MM$
is a finitely generated $\DD$-module. However, these assertions will be proved
directly in the following. In particular we will see that $\DD$ is a normal
semigroup ring over $K$.

By definition $\DD$ is a $\ZZ_+$-graded $R$-algebra with $\DD_0=M_0=R$, and
$\MM$ is a graded $\DD$-module if we assign degree $j$ to the elements of
$M_{jc+d}$. There exists $e>0$ such that $\DD$ is a finitely generated module
over its $R$-subalgebra generated by elements of degree $e$; for example, we
can take $e$ to be the least common multiple of the degrees of the generators
of $\DD$ as an $R$-algebra. Let $\EE$ be the $e$th Veronese subalgebra of
$\DD$. We decompose $\MM$ into the direct sum of its $\EE$-submodules
$$
\MM_k=\Dirsum_{j=0}^\infty M_{(ej+k)c+d},\qquad k=0,\dots,e-1.
$$
In view of what has to be proved, we can replace $\DD$ by $\EE$ and $\MM$ by
$\MM_k$. Then we have reached a situation in which $\DD$ is a finitely
generated module over the subalgebra generated by its degree $1$ elements.

Note that $\MM$ is isomorphic to an ideal of $\DD$: multiplication by a
monomial $X^a$ such that $a$ has residue class $-(d+k)$ in
$\Cl(R)\iso\ZZ^s/\sigma(\gp(S))$ maps $\MM$ into $\DD$. Since $\MM$ is not zero
(and $\DD$ is an integral domain), we see that $\Supp\MM=\Spec \DD$.

Let $\mm$ be the irrelevant maximal ideal of $R$; it is generated
by all elements $x\in S$, $x\neq1$ (in multiplicative notation).
Then clearly $\overline\MM=\MM/\mm\MM$ is a finitely generated
$\overline\DD=\DD/\mm\DD$-module. Note that $\overline\DD$ is a
$K$-algebra with $\overline\DD_0=K$ in a natural way. Furthermore
it is a finitely generated module over its subalgebra $\DDD$
generated by its degree $1$ elements. In particular $\overline\MM$
is a finitely generated $\DDD$-module. By construction (and
Nakayama's lemma) we have
$$
\mu(M_{(ej+k)c+d})=\dim_K M_{(ej+k)c+d}/\mm M_{(ej+k)c+d}=H(\overline\MM,j)
$$
where $H$ denotes the Hilbert function of $\overline\MM$ as a $\ZZ_+$-graded
$\overline\DD$- or $\DDD$-module. For $j\gg0$ the Hilbert function is
given by the Hilbert polynomial. It is a polynomial of degree $\delta-1$ where
$\delta$ is the Krull dimension of $\overline\MM$. Note that
$\Supp\overline\MM=\Spec \overline\DD$, since $\Supp\MM=\Spec \DD$; in
particular one has $\dim\overline\MM=\dim \overline\DD$. Moreover the leading
coefficient of the Hilbert polynomial is $e(\overline\MM)/(\delta-1)!$ and so
all the claims for $\MM$ follow if $m+1=\delta>1$.

At this point we have to clarify the structure of $\DD$ as a normal semigroup
ring over $K$. For convenience we choose a divisorial ideal $I\subset R$ of
class $c$ generated by monomials. Then there exists an $R$-module isomorphism
$M_{c}\to I$ mapping monomials to monomials, and such an isomorphism induces a
$K$-algebra isomorphism from $\DD$ to
$$
\R=\Dirsum_{j=0}^\infty I^{(j)}T^j\subset R[T]=K[S\dirsum \ZZ_+].
$$
There exist $a_1,\dots,a_s\ge0$ such that $I=\pp_1^{(a_1)}\cap
\dots \cap\pp_s^{(a_s)}$. The monomial corresponding to
$(u,z)\in\gp(S)\dirsum\ZZ$ belongs to $\R$ if and only if
$$
z\ge0,\qquad \sigma_i(u)-za_i\ge 0,\qquad i=1,\dots,s.
$$
It follows immediately that $\R$ is a normal semigroup ring over $K$. Let $\Sc$
be its semigroup of monomials. One has $\gp(\Sc)=\gp(S)\dirsum \ZZ$, and the
elements with last component $j$ give the monomials of $I^{(j)}$.

It is not hard to show that the faces of $C(\Sc)$ that are not
contained in $C(S)$ are the closed envelopes of the
$\RR_+$-envelopes of the faces of $C(I)'=\{(x,1)\: x\in C(I)\}$.

Moreover, exactly those faces $F$ that do not contain an element
from $\mm$ intersect $C(I)'$ in a compact face. In fact, if $F$
contains a monomial $x\in\mm$, then it contains $y+kx$, $k\in
\ZZ_+$, for each $y\in F$, and therefore an unbounded set. If $F$
does not contain an element of $\mm$, then the linear subspace
spanned by the elements of $S$ intersects $F$ in a single point,
and thus each translate intersects $F$ in a compact set.

Since the dimension of $\R/\mm\R$ is just the maximal dimension of
a face $F$ of $C(\Sc)$ not containing an element of $\mm$, we see
that $\dim\R/\mm\R=m+1$. In fact, the largest dimension of a
compact face of $C(I)'$ is $m$, and such a face extends to an
$m+1$-dimensional face of $C(\Sc)$.

Thus $\delta=m+1$ and $\delta>1$, since $C(I)$ has at least a $1$-dimensional
compact face: by hypothesis $I$ is not of torsion class.

For part (c) we note that $\htt\mm R=\grade \mm \R=\dim \R - \dim \R/\mm
\R=\dim R-m$ since $\R$ is Cohen--Macaulay by Hochster's theorem (and all the
invariants involved are stable under localization with respect to the maximal
ideal of $\R$ generated by monomials). Moreover $\grade \mm \R=\inf_j\depth
M_{cj}$, as follows by arguments analogous to those in the proof of Theorem
\ref{depth}.

By similar arguments the inequality for $\inf_j\depth M_{cj+d}$ results from
$\htt \mm\R=\dim R-m$.
\end{proof}

\begin{remark}\label{difflim}
The limits in Theorem \ref{aripro}(b) coincide if and only if the
$\DDD$-modules $\MM_k/\mm\MM_k$ all have the same multiplicity.
However, in general this is not the case. As an example one can take the
semigroup ring
$$
R=K[U^2,UV,V^2,XW,YW,XZ,YZ]\subset P=K[U,V,X,Y,Z,W]
$$
in its standard embedding. It has divisor class group $\Cl(R)=\ZZ/(2)\dirsum
\ZZ$. The non-zero torsion class is represented by the coset module
$M_{(1,0)}=RU+RV$, and $M_{(0,1)}=RX+RY$ represents a generator of the direct
summand $\ZZ$. Let $c\in\Cl(R)$ be the class of $M_{(1,1)}$. As an $R$-module,
$M_{jc}$, $j$ odd, is generated by the monomials $U\mu,V\mu$ where $\mu$ is a
degree $j$ monomial in $X,Y$, whereas for even $j$ the monomials $\mu$ form a
generating system. The limits for $k=0$ and $k=1$ therefore differ by a factor
of $2$ ($d=0$, $e=2$).
\end{remark}

\begin{proof}[Second proof of Theorem \ref{finite}]
Let $P$ be the polynomial ring of the standard embedding of $R$. Then $P$ is a
$\Cl(R)$-graded $R$-algebra whose graded component $P_c$, $c\in\Cl(R)$ is the
module $M_c$. Passing to residue classes modulo $\mm$ converts the assertion of
the theorem into a statement about the Hilbert function (with respect to $K$)
of the $\Cl(R)$-graded $K$-algebra $P/\mm P$; note that $(P/\mm P)_0=R/\mm=K$.
By Theorem \ref{aripro} the Hilbert function goes to infinity along each
arithmetic progression in $\Cl(R)$. Therefore we are in a position to apply
Theorem \ref{Hilbert} below. It says that there are only finitely many $c\in
\Cl(R)$ where $\mu(M_c)=H(P/\mm P,c)$ does not exceed a given bound $m$.
\end{proof}

This deduction of Theorem \ref{finite} uses the combinatorial hypotheses on $R$
only at a single point in the proof of Theorem \ref{aripro}, namely where we
show that $\dim\DD/\mm\DD \ge 2$. Thus the whole argument can be transferred
into a more general setting, provided an analogous condition on dimension
holds.

\section{Convex ideals}\label{gen}

Let $\xi_1,\dots,\xi_n$ be $\ZZ$-linear forms on $\ZZ^r$. Then
\begin{equation*}
S=\{x\in\ZZ^r\: \xi_i(x)\ge 0,\ i=1,\dots,n\}\tag{$*$}
\end{equation*}
is a normal semigroup. For $a\in\ZZ^n$ we set
$$
T(a;\xi)=\{x\in\ZZ^r\: \xi_i(x)\ge a_i,\ i=1,\dots,n\}.
$$
Let $K$ be a field. Then $I(a;\xi)=KT(a;\xi)$ is an $R=K[S]$-module in a
natural way.

In the following we will always assume that $S$ is positive and that
$\gp(S)=\ZZ^r$; furthermore we will assume that $\xi$ is non-degenerate, i.e.\
$\xi_i\neq0$ for all $i$; this assumption is certainly no restriction of
generality. Then it follows easily that $I(a;\xi)$ is a non-zero fractional
ideal of $R$ for every $a$. These ideals are called \emph{$\xi$-convex}.

If $\sigma=(\sigma_1,\dots,\sigma_s)$ consists of the support forms of $S$,
then the $\sigma$-convex ideals are just the monomial divisorial ideals of $R$
studied in the previous sections.

The linear map $\xi\: \gp(S)\to\ZZ^n$ evidently restricts to a pure embedding
of $S$ into $\ZZ^n$ (if $S$ is positive), and, conversely, a pure embedding
induces a presentation $(*)$ of $S$. As in the proof of Theorem \ref{Standard}
one sees that each coset modules of $\xi(S)$ is isomorphic to a $\xi$-convex
ideal, and vice versa. Proposition \ref{divi} therefore describes the condition
under which all $\xi$-convex ideals are divisorial (and in its proof we have
actually switched from a coset module to the corresponding $\xi$-convex ideal).

However, in contrast to the standard embedding, different coset modules can now
be isomorphic. Proposition \ref{iso} below clarifies this fact. For
$T\subset\gp(S)$ we define the ``effective $\xi$-bounds'' $\eff(T;\xi)\in\ZZ^n$
component wise by
$$
\eff(T;\xi)_i=\inf\{\xi_i(t)\: t\in T\};
$$
for a fractional monomial ideal $I$ we set $\eff(I;\xi)=\eff(T;\xi)$ where $T$
is the monomial basis of $I$; and for $a\in\ZZ^n$ we set
$\eff(a;\xi)=\eff(I(a;\xi);\xi)$. The convex polyhedron $C(I;\xi)\subset \RR^r$
associated with $I$ is given by
$$
C(I;\xi)=\{x\in \RR^r\: \xi_i(x)\ge \eff(I;\xi)_i,\ i=1,\dots,n\}.
$$

\begin{proposition}\label{iso}
\begin{itemize}
\item[(a)]
The assignment $I\mapsto C(I;\xi)$ is injective on the set of $\xi$-convex
ideals $I$.
\item[(b)]
The following are equivalent for $\xi$-convex ideals $I$ and $J$:
\begin{itemize}
\item[(i)] $I$ and $J$ are isomorphic $R$-modules.
\item[(ii)] There exists $y\in\gp(S)$ such that $\eff(J;\xi)=\eff(I,\xi)+\xi(y)$.
\item[(iii)] There exists $y\in\gp(S)$ such that $C(I;\xi)=C(J;\xi)+\xi(y)$.
\end{itemize}
\item[(c)] If $C(I;\xi)$ has a single extreme point for a $\xi$-convex ideal,
then $I$ is a divisorial ideal of torsion class.
\item[(d)]
The coset modules $Q_a,Q_b$ of $\xi(S)$ associated with the residue classes of
$a,b\in\ZZ^n$ are isomorphic if and only if $\eff(a;\xi)\equiv \eff(b;\xi)\mod
\xi(\gp(S))$.
\end{itemize}
\end{proposition}

\begin{proof}
(a) One can recover $I$ from $C(I;\xi)$.

(b) Fractional monomial ideals $I$ and $J$ of $K[S]$ are isomorphic if and only
if there exists a monomial $x\in K[\gp(S)]$ such that $J=Ix$. If one converts
this into additive notation, one immediately obtains the equivalence of (i),
(ii), and (iii).

(c) Set $a_i=\eff(I;\xi_i)$, $i=1,\dots,n$. Each of the hyperplanes given by
the equations $\xi_i(x)=a_i$ passes through an extreme point of $C(I;\xi)$.
Since there is a single extreme point $z$, all the parallels to the support
hyperplanes pass through this point. Thus $C(I;\xi)=C(S)+z$, and $C(I;\xi)$ is
cut out by parallels of the support hyperplanes of $S$. Then $I$ is a
divisorial ideal. Divisorial ideals $I$ for which $C(I)$ has a single extreme
point are of torsion class as observed in the proof of Theorem \ref{finite}.

(d) Since $Q_a\iso I(\eff(a;\xi);\xi)$ and $Q_b\iso I(\eff(b;\xi);\xi)$, this
follows from (b).
\end{proof}

It is now clear that Lemma \ref{length} generalizes to the case of $\xi$-convex
ideals, and since the isomorphism classes of $\xi$-convex ideals are again
parameterized by the equivalence classes of the polyhedra $C(I;\xi)$ modulo
translations by vectors from $\gp(S)$, we obtain a complete generalization of
Theorem \ref{finite}:

\begin{theorem}\label{finite2}
For every $m\in\ZZ_+$ there exist only finitely many isomorphism classes of
$\xi$-convex ideals such that a representative $I$ has $\mu(I)\le m$.
\end{theorem}

There is no need for a corollary regarding the Cohen--Macaulay property, since
Cohen--Macaulay fractional ideals are automatically divisorial. As outlined in
Remark \ref{finrem}, Theorem \ref{finite2} has a purely combinatorial statement
(and, by the way, has been proved by purely combinatorial arguments).

\begin{remark}\label{convrem}
(a) Let us set $\Eff(\ZZ^n;\xi)=\{\eff(a;\xi)\: a\in \ZZ^n\}$. Then it is not
hard to show that $\Eff(\ZZ^n;\xi)$ is a subsemigroup of $\ZZ^n$ containing
$\xi(\gp(S))$, and it follows immediately from Proposition \ref{iso} that the
semigroup $\Cl(R;\xi)=\Eff(\ZZ^n;\xi)/\xi(\gp(S))$ parametrizes the isomorphism
classes of $\xi$-convex ideals.

For the converse of Proposition \ref{iso}(c) one has to replace ``torsion
class'' (in $\Cl(R)$) by ``torsion class in $\Cl(R;\xi)$''.

Theorem \ref{aripro} generalizes to arithmetic progressions in $\Cl(R;\xi)$. We
leave the details to the reader. If $\Cl(R;\xi)$ is always finitely generated,
then the second proof of Theorem \ref{finite} can also be generalized.

(b) Theorem \ref{finite2} can be generalized to arbitrary normal affine
semigroups in the same way as Theorem \ref{finite} was generalized in Remark
\ref{finrem}.

(c) It is necessary to fix $\xi$ in Theorem \ref{finite2}. Already for
$S=\ZZ_+^2$ one can easily find infinitely many pairwise non-isomorphic ideals
generated by $2$ elements such that each of them is $\xi$-convex for a suitable
$\xi$.
\end{remark}

\section{On the growth of Hilbert functions}\label{Gr}

We introduce some terminology: if $S$ is a subsemigroup of an abelian group
$G$, then $T\subset G$ is an $S$-module if $S+T\subset T$ (the case
$T=\emptyset$ is not excluded). If $S$ is finitely generated and $T$ is a
finitely generated $S$-module, then every $S$-module $T'\subset T$ is also
finitely generated. For example, this follows by ``linearization'' with
coefficients in a field $K$: $M=KT\subset K[G]$ is a finitely generated module
over the noetherian ring $R=K[S]$, and so all its submodules are finitely
generated over $R$. For $KT'$ this implies the finite generation of $T'$ over
$S$.

First we note a result on the finite generation of certain subalgebras of
graded algebras and submodules of graded modules. We do not know of a reference
covering it in the generality of Theorem \ref{fingen}.

Gordan's lemma says that a subsemigroup $S\subset \ZZ^n$ is finitely generated
if and only if the cone $C(S)\subset\RR^n$ is finitely generated. It implies
that the integral closure $\hat S$ of an affine subsemigroup $S$ in $\ZZ^n$ is
finitely generated, since $C(S)=C(\hat S)$. The reader may check that we do not
use Gordan's lemma in the proof of Theorem \ref{fingen}. We will however use
that a finitely generated cone has only finitely many support hyperplanes and
that $\hat S=C(S)\cap \ZZ^n$.

\begin{theorem}\label{fingen}
Let $G$ be a finitely generated abelian group, $S$ a finitely generated
subsemigroup of $G$, and $T\subset G$ a finitely generated $S$-module.
Furthermore let $R$ be a noetherian $G$-graded ring and $M$ a $G$-graded
finitely generated $R$-module. Then the following hold:
\begin{itemize}
\item[(a)]
$R_0$ is noetherian ring, and each graded component $M_g$, $g\in G$, of $M$ is
a finitely generated $R_0$--module.
\item[(b)] $A=\Dirsum_{s\in S}R_s$ is a finitely generated $R_0$-algebra.
\item[(c)] $N=\Dirsum_{t\in T} M_t$ is a finitely generated $A$-module.
\end{itemize}
\end{theorem}

\begin{proof}
(a) One easily checks that $M'R\cap M_g=M'$ for each $R_0$-submodule $M'$ of
$M_g$. Therefore ascending chains of such submodules $M'$ of $M_g$ are
stationery.

(b) First we do the case in which $G$ is torsionfree, $G=\ZZ^m$, and $S$ is an
integrally closed subsemigroup of $\ZZ^m$.

Let $\phi\:\ZZ^m\to \ZZ$ be a non-zero linear form. It induces a $\ZZ$-grading
on $R$ with $\deg_{\ZZ}(a)=\phi(\deg_{\ZZ^M}(a))$ for each non-zero
$\ZZ^m$-homogeneous element of $R$. Let $R'$ denote $R$ with this
$\ZZ$-grading. Set $R'_-=\Dirsum_{k\le0}R'_k$ and define $R'_+$ is analogously.
By \cite[1.5.5]{BH} the $R'_0$-algebras $R'_+$ and $R'_-$ are finitely
generated $R'_0$-algebras, and $R'_0$ is a noetherian ring (by (a)). On the
other hand, $R'_0$ is a $(\Ker\phi)$-graded ring in a natural way, and by
induction we can conclude that $R'_0$ is a finitely generated $R_0$-algebra.

If $S=\ZZ^m$, then it follows immediately that $R$, the sum of $R_-$ and $R_+$
as an $R_0$-algebra, is again a finitely generated $R_0$-algebra.

Otherwise $S=\ZZ^m\cap C(S)$, and $C(S)$ has at least one support hyperplane:
$$
S=\{s\in \ZZ^m\: \alpha_i(s)\ge 0,\ i=1,\dots,v\}
$$
with $v\ge 1$. We use induction on $v$, and the induction hypothesis applies to
$R'=\Dirsum_{s\in S'}R_s$,
$$
S'= \{s\in \ZZ^m\: \alpha_i(s)\ge 0,\ i=1,\dots,v-1\}.
$$
Applying the argument above with $\phi=\alpha_v$, one concludes that $A=R'_+$
is a finitely generated $R_0$-algebra.

In the general case for $G$ and $S$ we set $G'=G/H$ where $H$ is the torsion
subgroup of $G$, and denote the natural surjection by $\pi\:G\to G'$. Let $R'$
be $R$ with the $G'$-grading induced by $\pi$ (its homogeneous components are
the direct sums of the components $R_g$ where $g$ is in a fixed fiber of
$\pi$). Let $S'$ be the integral closure of $\pi(S)$ in $G'$. Then
$A'=\Dirsum_{s'\in S'}R'_{s'}$ is a finitely generated algebra over the
noetherian ring $R'_0$, as we have already shown. But $R'_0$ is a finitely
generated module over $R_0$ by (a), and so $A'$ is a finitely generated
$R_0$-algebra. In particular, $R$ itself is finitely generated over $R_0$.

It is not hard to check that $A'$ is integral over $A$; in fact, each element
$s\in \pi^{-1}(S')$ has a power $s^n\in S$ for suitable $n\in \NN$. Furthermore
it is a finitely generated $A$-algebra, and so a finitely generated $A$-module.
But then a lemma of Artin and Tate (see Eisenbud \cite[p.\ 143]{Ei}) implies
that $A$ is noetherian. As shown above, noetherian $G$-graded rings are
finitely generated $R_0$-algebras.

(c) By hypothesis, $T$ is the union of finitely many translates $S+t$.
Therefore we can assume that $T=S+t$. Passing to the shifted module $M(-t)$
(given by $M(-t)_g=M_{g-t}$), we can even assume that $S=T$. Now the proof
follows the same pattern as that of (b). In order to deal with an integrally
closed subsemigroup of a free abelian group $G=\ZZ^m$, one notes that $M_+$ is
a finitely generated module over $R_+$ where $M_+$ is the positive part of $M$
with respect to a $\ZZ$-grading (induced by a linear form $\phi\:\ZZ^m\to\ZZ$).
This is shown as follows: the extended module $RM_+$ is finitely generated over
$R$ , and every of its generating systems $E\subset M_+$ together with finitely
many components $M_i$, $i\ge0$, generate $M_+$ over $R_+$; furthermore the
$M_i$ are finitely generated over $R_0$ by (a).

For the general situation we consider $N'$ defined analogously as $A'$. It is a
finitely generated $A'$-module by the previous argument. Since $A'$ is a
finitely generated $A$-module, $N'$ is finitely generated over $A$, and so is
its submodule $N$.
\end{proof}

We note a purely combinatorial consequence.

\begin{corollary}\label{fincomb}
Let $S$ and $S'$ be affine subsemigroups of $\ZZ^m$, $T\subset
\ZZ^m$ a finitely generated $S$-module, and $T'\subset\ZZ^m$ a
finitely generated $S'$-module. Then $S\cap S'$ is an affine
semigroup, and $T\cap T'$ is a finitely generated $S\cap
S'$-module.
\end{corollary}

\begin{proof}
We choose a field $K$ of coefficients and set $R=K[S']$, $M=KT'$. Then the
hypotheses of the theorem are satisfied, and it therefore implies the finite
generation of
$$
A=\Dirsum_{s\in S} R_s=K[S\cap S'],\qquad\text{and}\qquad N=\Dirsum_{t\in
T}M_t=K(T \cap T').
$$
as a $K=R_0$-algebra and an $A$-module respectively. However, finite generation
of the ``linearized'' objects is equivalent to that of the combinatorial ones.
\end{proof}

The next theorem is our main result on the growth of Hilbert functions. Note
that we do not assume that $R_0=K$; the graded components of $R$ and $M$ may
even have infinite $K$-dimension.

\begin{theorem}\label{Hilbert}
Let $K$ be a field, $G$ a finitely generated abelian group, $R$ a
noetherian $G$-graded $K$-algebra for which $R_0$ is a finitely
generated $K$-algebra, and $M$ a finitely generated $G$-graded
$R$-module. Consider a finitely generated subsemigroup $S$ of $G$
containing the elements $\deg r$, $r\in R\setminus\{0\}$
homogeneous, and a finitely generated $S$-submodule $T$ of $G$
containing the elements $\deg x$, $x\in M\setminus\{0\}$
homogeneous. Furthermore let $H$ be the $G$-graded Hilbert
function, $H(M,t)=\dim_K M_t$ for all $t\in G$.

Suppose $\lim_{k\to\infty} H(M,kc+d)=\infty$ for all choices of $c\in S$, $c$
not a torsion element of $G$, and $d\in T$. Then
$$
\#\{t\in T: H(M,t)\le C\} <\infty
$$
for all $C\in\ZZ_+$.
\end{theorem}

Note that $R$ is a finitely generated $R_0$-algebra by Theorem
\ref{fingen}, and therefore a finitely generated $K$-algebra. Let
$S'$ be the subsemigroup of $G$ generated by the elements $\deg
r$, $r\in R\setminus\{0\}$ homogeneous, and $T'$ be the
$S'$-submodule of $G$ generated by the elements $\deg x$, $x\in
M\setminus\{0\}$ homogeneous. Then all the hypotheses are
satisfied with $S'$ in place of $S$ and $T'$ in place of $T$.
However, for technical reasons the hypothesis of the theorem has
to be kept more general. (We are grateful to Robert Koch for
pointing out some inaccuracies in previous versions of the theorem
and its proof.)

\begin{proof}[Proof of Theorem \ref{Hilbert}]
We split $G$ as a direct sum of a torsionfree subgroup $L$ and its
torsion subgroup $G_{\mathrm{tor}}$. Let $R'=\Dirsum_{\ell\in L}
R_\ell$, and split $M$ into the direct sum
$$
M=\Dirsum_{h\in G_{\mathrm{tor}}} M'_h,\qquad M'_h=\Dirsum_{\ell\in L}
M_{(\ell,h)}.
$$
By Theorem \ref{fingen} $R'$ is a finitely generated $R_0$-algebra and $M'_h$
is a finitely generated $L$-graded $R'$-module for all $h\in G_{\mathrm{tor}}$,
and since the hypothesis on the Hilbert function is inherited by $M'_h$, it is
enough to do the case $G=L=\ZZ^m$.

We use induction on $m$. In the case $m=1$ it is not difficult to
see (and well-known) that $T$ is the union of finitely many
arithmetic progressions that appear in the hypothesis of the
theorem.

As a first step we want to improve the hypothesis on the Hilbert function from
a ``$1$-dimensional '' condition to a ``$1$-codimensional'' condition by an
application of the induction hypothesis.

Let $U$ be a subgroup of $L$ with $\rank U<\rank L$ and $u\in L$. Then $U$ is finitely
generated as a subsemigroup. We set
$$
R'=\Dirsum_{s\in U} R_s\quad\text{and}\quad M'=\Dirsum_{t\in U+u} M_t.
$$
Theorem \ref{fingen} implies that $R'$ is a finitely generated $K$-algebra, and
$M'$ is a finitely generated $R'$-module.

After fixing an origin in $U+u$ we can identify it with $U$. Therefore we can
apply the induction hypothesis to $R'$ and $M'$. It follows that
\begin{equation*}
\#\{t\in T\cap (U+u): H(t,M)\le C\} <\infty.\tag{$*$}
\end{equation*}

By Theorem \ref{fingen} $R$ is a finitely generated $R_0$-algebra
and thus a finitely generated $K$-algebra. We represent $R$ as the
residue class ring of an $L$-graded polynomial ring $P$ over $K$
in a natural way (in particular the monomials in $P$ are
homogeneous in the $L$-grading). The hypothesis that $R_0$ is a
finitely generated $K$-algebra is inherited by $P$ since $P_0$ is
a (not necessarily positive) normal affine semigroup ring. Thus we
may assume that $R$ itself is generated by finitely many
algebraically independent elements as a $K$-algebra.

Obviously $M$ has a filtration
$$
0=M_0\subset M_1\subset\dots\subset M_n=M
$$
where each successive quotient $M_{i+1}/M_i$ is a cyclic $L$-graded $R$-module,
that is
$$
M_{i+1}/M_i\iso (R/I_{i+1})(-s_i)
$$
with an $L$-graded ideal $I_{i+1}$ in $R$ and a shift $s_i\in L$. As far as the
Hilbert function is concerned, we can replace $M$ by the direct sum of these
cyclic modules. After the introduction of a term order we can replace
$R/(I_{i+1})(-s_i)$ by $R/(\ini(I_{i+1}))(-s_i)$ where $\ini(I_{i+1})$ is the
initial ideal. It is well known that $R/\ini(I_{i+1})$ has a filtration whose
successive quotients are of the form $R/\pp$ with a prime ideal $\pp$ generated
by monomials, and therefore by indeterminates of $R$. (For example, see the
proof of \cite[4.1.3]{BH}, and use that associated prime ideals of multigraded
modules are generated by indeterminates if the multigrading is that induced by
the semigroup of all monomials in $R$.)

Altogether this reduces the problem to the case in which the $K$-vector space
$M$ is isomorphic to the direct sum of vector spaces $P_i(-s_i)$ where $P_i$ a
polynomial ring generated by indeterminates with degrees in $L$, and $s_i\in
L$. Furthermore we can use that the Hilbert function of $M$ satisfies condition
$(*)$. The Hilbert function now counts the total number of monomials in each
degree. Replacing the monomials by their exponent vectors, we can deduce the
theorem from the next one.
\end{proof}

\begin{theorem}\label{Hilbcomb}
Let $G$ be a finitely generated abelian group, $S$ a finitely generated
subsemigroup of $G$, and $T$ a finitely generated $S$-submodule of $G$.
Consider maps
$$
\psi_i\:A_i\to T,\qquad \psi_i(x)=\phi_i(x)+t_i\ \text{for all $x\in A_i$}
$$
where $A_i$ is an affine semigroup, $\phi_i:A_i\to S$ is a homomorphism of
semigroups, and $t_i\in T$, $i=1,\dots,v$. Furthermore let
$$
\Psi\:A_1\cup\dots\cup A_v\to T,\qquad \Psi|_{A_i}=\psi_i,
$$
be the map defined on the disjoint union of the $A_i$ by all the $\psi_i$. For
$t\in G$ set
$$
H(t)=\#\{x\in A_1\cup\dots\cup A_v:\Psi(x)=t\}.
$$
Suppose that $\lim_{k\to\infty} H(kc+d)=\infty$ for all $c\in S$, $c$ not a
torsion element of $G$, and all $d\in T$. Then
$$
\#\{t\in T: H(t)\le C\} <\infty
$$
for all $C\in\ZZ_+$.
\end{theorem}

\begin{proof}
In step (a) we prove the theorem under the assumption that $G=L=\ZZ^m$ for some
$m$ and that
$$
\#\{t\in T\cap (U+u): H(t)\le C\} <\infty.
$$
for all subgroups $U$ of $L$ with $\rank U<\rank L$. This is enough to complete the proof of
Theorem \ref{Hilbert}. In step (b) we can then use Theorem \ref{Hilbert}.

(a) The first observation is that we can omit all the maps $\psi_i$ that are
injective. This reduces the function $H$ in each degree by at most $v$, and has
therefore no influence on the hypothesis or the desired conclusion.

The difficult case is $C=0$, and we postpone it. So suppose that we have
already shown that the number of ``gaps'' (elements in $T$ with no preimage at
all) is finite. Then we can restrict ourselves to $\Im \Psi$ if we want to show
that there are only finitely many elements with at most $C>0$ preimages.

It is enough to show that the elements in $\Im \psi_i$ with at most $C$
preimages are contained in the union of finitely many sets of the form $U+u$
where $U$ is a proper direct summand of $L$. Then we can use the hypothesis on
the sets $\{x\in T\cap (U+u): H(x)\le C\}$. We can certainly assume that $v=1$
and $t_1=0$, and have only to consider a \emph{non-injective, surjective}
homomorphism $\phi:A\to S$.

For an ideal (i.~e.\ $S$-submodule) $I\neq\emptyset$ of $S$ we have that
$S\setminus I$ is contained in finitely many sets $U+u$. In fact, $I$ contains
an ideal $J\neq\emptyset$ of the normalization $\bar S$ of $S$, namely the
conductor ideal $\{s\in S\: \bar S+s\subset S\}$. (This follows from the
corresponding theorem of commutative algebra; see \cite[Lemma 5.3]{Gu2} where
the assertion has been proved under the superfluous condition that $S$ is
positive.) Therefore $S$ contains a set $\bar S+s$ with $s\in S$. It follows
that $\bigl(S\setminus I\bigr)\subset \bigl(\bar S\setminus(\bar S+s)\bigr)$.
The latter set is contained in finitely many parallels to the facets of $S$
through points of $S$ (in $\gp(S)$). Each of these parallels is itself
contained in a set of type $U+u$. To sum up, it is enough to find an ideal $I$
in $S$ such that each element of $I$ has at least $C+1$ preimages.

Now we go to $A$ and choose $a\in A$ such that $\bar A+a\subset A$ where $\bar
A$ is again the normalization. The homomorphism $\phi$ has a unique extension
to a group homomorphism $\gp(A)\to L$, also denoted by $\phi$. By assumption
$\Ker\phi\neq 0$. A sufficiently large ball $B$ in $\gp(A)\tensor \RR$ with
center $0$ therefore contains $C+1$ elements of $\Ker\phi$, and there exists
$b\in \bar A$ for which $B+b$ is contained in the cone $\RR_+ A$. Thus
$(B\cap\Ker\phi)+b\subset \bar A$. It follows that each element in $I=\phi(\bar
A+a+b)$ has at least $C+1$ preimages. Since $\bar A+a+b$ is an ideal in $A$ and
$\phi$ is surjective, $I$ is an ideal in $S$.

(b) By linearization we now derive Theorem \ref{Hilbcomb} from Theorem
\ref{Hilbert}. Let $K$ be a field. Then we set $R_i=K[A_i]$, and the
homomorphism $\phi_i$ allows us to consider $R_i$ as a $G$-graded $K$-algebra.
Next we choose a polynomial ring $P_i$ whose indeterminates are mapped to a
finite monomial system of generators of $A_i$, and so $P_i$ is also $G$-graded.
Set
$$
R=P_1\tensor_K\dots\tensor_K P_v\qquad\text{and}\qquad
M=R_1(-s_1)\dirsum\dots\dirsum R_v(-s_v)
$$
Evidently $R$ is a finitely generated $G$-graded $K$-algebra; in particular it
is noetherian. Moreover $R_i$ is residue class ring of $R$ in a natural way,
and therefore $R_i(-s_i)$ can be considered a $G$-graded $R$-module. Therefore
$M$ is a $G$-graded $R$-module whose Hilbert function is the function $H$ of
the theorem.
\end{proof}

It remains to do the case $C=0$. For simplicity we only formulate it under the
special assumptions of step (a) in the proof of Theorem \ref{Hilbcomb}. We
leave the general as well as the commutative algebra version to the reader. The
semigroups $A_i$ of Theorem \ref{Hilbcomb} can now be replaced by their images.

\begin{proposition}\label{cover}
Let $L=\ZZ^m$, $S$ an affine subsemigroup of $L$, $T$ a finitely generated
$S$-submodule of $L$. Consider subsemigroups $A_1,\dots,A_v$ of $L$ and
elements $t_1,\dots,t_v\in T$ such that the set
$$
\Gc=T\setminus ((A_1+t_1)\cup\dots\cup(A_v+t_v))
$$
of ``gaps'' satisfies the following condition: for each subgroup $U$ of
$L$ with $\rank U<\rank L$ and each $u\in L$ the intersection $(U+u)\cap \Gc$ is finite. Then $\Gc$ is
finite.
\end{proposition}

\begin{proof}
Note that $T$ is contained in finitely many residue classes modulo $\gp(S)$.
Therefore we can replace each $A_i$ by $A_i\cap\gp(S)$: the intersection of
$A_i+t_i$ with a residue class modulo $\gp(S)$ is a finitely generated
$A_i\cap\gp(S)$-module by Corollary \ref{fincomb}.

We order the $A_i$ in such a way that $A_1,\dots,A_w$ have the same rank as
$S$, and $A_{w+1},\allowbreak\dots,\allowbreak A_v$ have lower rank. Let $W$ be
the intersection of $\gp(A_i)$, $i=1,\dots,w$. Since $\gp(S)/W$ is a finite
group, we can replace all the semigroups involved by their intersections with
$W$, and split the modules into their intersection with the residue classes
modulo $W$. We have now reached a situation where $A_i\subset\gp(S)$ for all
$i$, and $\gp(A_i)=\gp(S)$, unless $\rank A_i< \rank S$.

Next one can replace the $A_1,\dots,A_w$ by their normalizations. In this way
we fill the gaps in only finitely many $U+u$ (compare the argument in the proof
of Theorem \ref{Hilbcomb}), and therefore we fill only finitely many gaps.

At this point we can assume that $A_1,\dots,A_w$ are integrally closed in $L$.
Furthermore we must have $C(S)\subset C(A_1)\cup\dots\cup C(A_w)$ -- otherwise
an open subcone of $C(S)$ would remain uncovered, and this would remain so in
$T$: the lower rank semigroups cannot fill it, and neither can it be filled by
finitely many translates $U+u$ where $U$ is a subsemigroup of $S$ with $\rank U<\rank S$. In
fact, $(A_i+t_i)\setminus A_i$ is contained in the union of finitely many such
translates, and the same holds for $(C(S)\cap L)\setminus S$. Since
$A_1,\dots,A_w$ are integrally closed, we have $S\subset A_1\cup\dots\cup A_w$.

Now we choose a system of generators $u_1,\dots,u_q$ of $T$ over $S$. We have
$$
T\subset\bigcup_{i,j} A_i+u_j.
$$
But $A_i+u_j$ and $A_i+t_i$ only differ in finitely many translates of proper
direct summands of $L$ parallel to the support hyperplanes of $A_i$. So in the
last step we have filled only finitely many gaps. Since no gaps remain, their
number must have been finite from the beginning.
\end{proof}

\end{document}